\documentclass[12pt,oneside]{amsart}

\usepackage{amsfonts}
\usepackage{amscd}
\usepackage{amsthm}
\usepackage{amssymb}
\usepackage{amsmath}
\usepackage{epsfig}

\usepackage{graphicx}
\usepackage{xypic}
\input xypic
\input xy
\xyoption{all}

\newcommand{\mD}{\mathbb D}
\newcommand{\mR}{\mathbb R}

\newcommand{\tX}{\tilde X}
\newcommand{\pinf}{\partial ^\infty}
\newcommand{\p}{\prime}

\newcommand{\dlim}{\underrightarrow{\lim}}

\title{A boundary version of Cartan-Hadamard and applications to rigidity.}

\author{Jean-Fran\c cois Lafont}
\address{Department of Mathematics,
The Ohio State University, Columbus, OH 43210}
\email{jlafont@math.ohio-state.edu}

\theoremstyle{definition}
\newtheorem{Def}{Definition}[section]

\theoremstyle{plain}

\newtheorem{Thm}{Theorem}[section]
\newtheorem{Cor}{Corollary}[section]
\newtheorem{Lem}{Lemma}[section]
\newtheorem{Prop}{Proposition}[section]

\theoremstyle{remark}

\newtheorem*{Prf}{Proof}
\newtheorem*{Rmk}{Remark}

\begin{document}

\begin{abstract}
In this paper, we prove a version of the classical Cartan-Hadamard
theorem for negatively curved manifolds, of dimension $n\neq 5$,
with non-empty totally geodesic boundary.  More precisely, if
$M_1^n,M_2^n$ are any two such manifolds, we show that (1) $\pinf
\tilde M_1^n$ is homeomorphic to $\pinf \tilde M_2^n$, and (2)
$\tilde M_1^n$ is homeomorphic to $\tilde M_2^n$.
As a sample application, we show that simple, thick, negatively
curved P-manifolds of dimension $\geq 6$ are topologically rigid. We
include some straightforward consequences of topological rigidity
(diagram rigidity, weak co-Hopf property, and Nielson realization
problem).
\end{abstract}

\maketitle

\section{Introduction.}

A key aspect in the study of non-positively curved Riemannian
manifolds is the large number of {\it rigidity theorems} known to
hold for these spaces.  Two outstanding such theorems are (1) Mostow
rigidity \cite{M}, stating that in dimension $\geq 3$, homotopy
equivalence of irreducible locally symmetric spaces of non-compact
type implies isometry of the spaces, and (2) Farrell-Jones
topological rigidity \cite{FJ}, stating that in dimension $\geq 5$,
homotopy equivalence of non-positively curved Riemannian manifolds
implies homeomorphism of the spaces.

A natural question is how to extend these theorems to the context of
singular spaces satisfying a metric analogue of ``non-positive
curvature''.  In a series of papers (\cite{L1}, \cite{L2},
\cite{L3}), the author introduced the class of hyperbolic
P-manifolds, which one can view as some of the simplest non-manifold
CAT(-1) spaces, and established Mostow rigidity within this class of
spaces.  In the present paper, our interest lies in establishing
topological rigidity for negatively curved P-manifolds. Let us first
recall the definition of a P-manifold:

\begin{Def}
A closed $n$-dimensional {\it piecewise manifold} (henceforth
abbreviated to P-manifold) is a topological space which has a
natural stratification into pieces which are manifolds. More
precisely, we define a $1$-dimensional P-manifold to be a finite
graph. An $n$-dimensional P-manifold ($n\geq 2$) is defined
inductively as a closed pair $X_{n-1}\subset X_n$ satisfying the
following conditions:

\begin{itemize}
\item Each connected component of $X_{n-1}$ is either an
$(n-1)$-dimensional P-manifold, or an $(n-1)$-dimensional manifold.
\item The closure of each connected component of $X_n-X_{n-1}$
is homeomorphic to a compact orientable $n$-manifold with boundary,
and the homeomorphism takes the component of $X_n-X_{n-1}$ to the
interior of the $n$-manifold with boundary; the closure of such a
component will be called a {\it chamber}.
\end{itemize}

\noindent Denoting the closures of the connected components of
$X_n-X_{n-1}$ by $W_i$, we observe that we have a natural map $\rho:
\coprod
\partial W_i\longrightarrow X_{n-1}$ from the disjoint union
of the boundary components of the chambers to the subspace
$X_{n-1}$.  We also require this map to be surjective, and a
homeomorphism when restricted to each component. The P-manifold is
said to be \emph{thick} provided that each point in $X_{n-1}$ has at
least three pre-images under $\rho$. We will henceforth use a
superscript $X^n$ to refer to an $n$-dimensional P-manifold, and
will reserve the use of subscripts $X_{n-1},\ldots ,X_1$ to refer to
the lower dimensional strata.  For a thick $n$-dimensional
P-manifold, we will call the $X_{n-1}$ strata the {\it branching
locus} of the P-manifold.
\end{Def}

Intuitively, we can think of P-manifolds as being ``built'' by
gluing manifolds with boundary together along lower dimensional
pieces. Examples of P-manifolds include finite graphs and soap
bubble clusters.  Observe that compact manifolds can also be viewed
as (non-thick) P-manifolds. Less trivial examples can be constructed
more or less arbitrarily by finding families of manifolds with
homeomorphic boundary and glueing them together along the boundary
using arbitrary homeomorphisms. We now define the family of metrics
we are interested in.

\begin{Def}
A Riemannian metric on a 1-dimensional P-manifold (finite graph) is
merely a length function on the edge set.  A Riemannian metric on an
$n$-dimensional P-manifold $X^n$ is obtained by first building a
Riemannian metric on the $X_{n-1}$ subspace, then picking, for each
$W_i$ a Riemannian metric with non-empty totally geodesic boundary
satisfying that the gluing map $\rho$ is an isometry. We say that a
Riemannian metric on a P-manifold is {\it negatively curved} if at
each step, the metric on each $W_i$ is negatively curved.
\end{Def}

Observe that, at the cost of scaling the metric of the P-manifold
$X$ by a constant, one can assume that the metric on each $W_i$
has sectional curvature bounded above by $-1$.  Such a metric on the
P-manifold will automatically be locally CAT(-1), and hence the
fundamental group of a negatively curved P-manifold is a
$\delta$-hyperbolic group.  In particular, the universal cover
$\tilde X$ has a well-defined boundary at infinity, denoted $\pinf
\tilde X$.

\begin{Def}
We say that an $n$-dimensional P-manifold $X^n$ is {\it simple}
provided its codimension two strata is empty.  In other words, the
$(n-1)$-dimensional strata $X_{n-1}$ consists of a disjoint union of
$(n-1)$-dimensional manifolds.  We further assume that, for each
chamber $W_i$, the various boundary components of $W_i$ get attached
to distinct components of the codimension one strata.
\end{Def}

We can now state our main result:

\begin{Thm}[Topological rigidity]
Let $X_1,X_2$ be a pair of simple, thick, negatively curved
P-manifolds, of dimension $\geq 6$.  If $\pi_1(X_1)$ is isomorphic
to $\pi_1(X_2)$, then $X_1$ is homeomorphic to $X_2$.
\end{Thm}

We note that, corresponding to the stratification of a negatively
curved P-manifold, there is a natural diagram of groups having the
property that the direct limit of the diagram is precisely the
fundamental group of the P-manifold (by the generalized Seifert-Van
Kampen theorem).  Immediate consequences of the topological rigidity
are the following:

\begin{Cor}[Diagram rigidity]
Let $\mathcal D_1, \mathcal D_2$ be a pair of diagrams of groups,
corresponding to a pair of negatively curved, simple, thick
P-manifolds of dimension $n\geq 6$.  Then $\dlim \mathcal D_1$ is
isomorphic to $\dlim \mathcal D_2$ if and only if the two diagrams
are isomorphic.
\end{Cor}

\begin{Cor}[weak Co-Hopf property]
Let $X$ be a simple, thick, negatively curved P-manifold of
dimension $n\geq 6$, and assume that at least one of the chambers
has a non-zero characteristic number.  Then $\Gamma=\pi_1(X)$ is
weakly co-Hopfian, i.e. every injection $\Gamma\hookrightarrow \Gamma$
with image of finite index is in fact an isomorphism.
\end{Cor}

\begin{Cor}[Nielson realization problem]
Let $X$ be a simple, thick, negatively curved P-manifold of
dimension $n\geq 6$, and $\Gamma=\pi_1(X)$.  Then the
canonical map $Homeo(X)\rightarrow Out(\Gamma)$ is
surjective.
\end{Cor}

Now recall that a consequence of the classical Cartan-Hadamard
theorem is that if $M_1,M_2$ are a pair of closed $n$-dimensional
manifolds of non-positive sectional curvature, then the universal
covers $\tilde M_1$ and $\tilde M_2$ are homeomorphic (indeed, are
both diffeomorphic to $\mathbb R^n$). Another classic result is that
for such a manifold $M$, the boundary at infinity of the universal
cover $\pinf \tilde M$ is always homeomorphic to an $(n-1)$-sphere
$S^{n-1}$. The key to the proof of the previous two theorems is the
following analogue of these classic results, in the setting where
one allows a non-empty, totally geodesic boundary.

\begin{Thm}[Cartan-Hadamard]
Assume $M_1$, $M_2$ are a pair of compact, negatively curved
Riemannian manifolds of dimension $n\neq 5$, with non-empty, totally
geodesic boundary.  Then we have:
\begin{enumerate}
\item $\pinf \tilde M_1$ is homeomorphic to $\pinf \tilde M_2$.
\item $\tilde M_1$ is homeomorphic to $\tilde M_2$.
\end{enumerate}
where $\tilde M_i$ is the universal cover of $M_i$.
\end{Thm}

Note that if $n=2$, then the boundaries at infinity of the $\tilde
M_i$ are Cantor sets, and the first statement in the Theorem is just
the classical fact that any two Cantor sets are homeomorphic
(Brouwer's characterization theorem).

We now outline the organization of this paper.  In Section 2, we
will give a proof of Theorem 1.2.  The argument relies heavily on a
characterization of $n$-dimensional Sierpinski curves ($n\neq 3$)
due to Cannon \cite{C}.  The dimension restriction in Theorem 1.2
arises from the corresponding dimension restriction in Cannon's
work.

In Section 3, we will give outlines of the proofs of Theorems 1.1,
as well as proofs of the three corollaries.  The arguments for these follow
almost verbatim from
previous results of the author \cite{L2}, \cite{L3}.  More
precisely, in \cite{L2} the author gave a {\it topological} argument
allowing, in the case where the simple, thick, negatively curved
P-manifold was actually {\it hyperbolic} (i.e. each chamber is
isometric to a hyperbolic manifold with non-empty, totally geodesic
boundary), for recognizing the fundamental groups of the various
chambers and how they are attached together.  The argument in
\cite{L2} made use of the topology of the boundary at infinity of a
{\it hyperbolic} manifold with non-empty totally geodesic boundary.
Statement (1) of Theorem 1.2, and the fact that the argument in
\cite{L2} is purely topological, allows the entire argument to be
transferred to the case of simple, thick, {\it negatively curved}
P-manifolds, of dimension $n\neq 5$.

Once one knows how to recognize the fundamental groups of the
chambers and how they attach together, we can appeal, in dimension
$n\geq 5$, to the celebrated {\it Topological Rigidity Theorem} of
Farrell and Jones \cite{FJ}.  This allows us to completely determine
the topology (up to homeomorphism) of all the chambers, as well as
the topology of the codimension one strata.  Putting this together
yields the desired Theorem 1.1, and translating the topological
rigidity result (via the generalized Seifert-Van Kampen) into group
theoretic language immediately gives Corollary 1.1.
Exploiting the
correspondence between subgroups of fundamental groups and coverings
of the corresponding space, it is easy to obtain the weak co-Hopf
property (Corollary 1.2).  The Nielson realization problem
(Corollary 1.3) is immediate from the Theorem 1.1.

\begin{Rmk}
(1)  We note that topological rigidity fails (trivially) in
dimension $n=1$.  In dimension $n=2$, topological rigidity was
proved in \cite{L3}.  In dimension $n=3$, the argument given in the
present paper could be extended, provided one had an analogue of
Farrell-Jones \cite{FJ} for $3$-dimensional manifolds. This analogue
is a well-known consequence of the Geometrization Conjecture of
Thurston. A proof of the Geometrization Conjecture was announced a
few years ago by G. Perelman.  The author does not know whether
topological rigidity is to be expected in dimensions $n=4,5$ (though
the failure of the present proof is due to different reasons in each
of these two cases).

(2)  It would be interesting to see whether, in statement (2) of
Theorem 1.2, one can replace ``homeomorphism'' by
``diffeomorphism''.  As the reader will see, the argument in the
present paper has no chance of being extended to yield smooth
rigidity.

(3)  Concerning the hypothesis in Corollary 1.2 on the existence of a
non-zero characteristic number for one of the chambers, we point out
that the famous {\it Hopf Conjecture} on the sign of the Euler
characteristic asserts that for a closed, negatively curved, even
dimensional manifold $M^{2n}$, we have the inequality
$(-1)^n\chi(M^{2n})>0$.  It is easy to see (using a doubling
argument) that the Hopf Conjecture, if true, implies that for any
compact negatively curved manifold $M$ with {\it non-empty} totally
geodesic boundary, we have $\chi(M)\neq 0$.  In particular, the
validity of the Hopf Conjecture would yield the desired non-zero
characteristic number.  We also point out that a much stronger
result is known, namely Sela \cite{Se} has shown that a
non-elementary $\delta$-hyperbolic group is co-Hopfian if and only if
if is freely indecomposable.
\end{Rmk}

\vskip 5pt

\centerline {\bf Acknowledgements}

\vskip 5pt

The author would like to thank M. Hindawi and T. Januszkiewicz for
raising the question of determining the topology of the universal
cover of a compact non-positively curved Riemannian manifold with
non-empty totally geodesic boundary.  Theorem 1.2(b) provides a
partial answer to their question, in the case where the dimension
$n\neq 5$ and where one has the stronger assumption that the
sectional curvatures are strictly negative (though see the remarks
at the end of Section 2).

\section{Cartan-Hadamard for manifolds with boundary.}

We now proceed to prove Theorem 1.2 from the introduction.  So let
$M_1,M_2$ be a pair of compact, negatively curved manifolds of
dimension $n\neq 5$, with non-empty totally geodesic boundary.  We
will start by establishing property (1), namely that $\pinf \tilde
M_1$ is homeomorphic to $\pinf \tilde M_2$.  In order to do this, we
will make use of the characterization of Sierpinski curves due to
Cannon \cite{C} (generalizing a classic result of Whyburn \cite{W} in
dimension $n=2$).  We first start with a definition:

\begin{Def}
Let $\{U_i\}$ be a countable collection of pairwise disjoint subsets
of $S^n$ satisfying the following four conditions:
\begin{enumerate}
\item the collection $\{U_i\}$ forms a null sequence, i.e. $\lim \{diam(U_i)\}=0$,
\item $S^n-U_i$ is an $n$-cell for each $i$,
\item $Cl(U_i)\cap Cl(U_j)=\emptyset$ for each $i\neq j$ ($Cl$ denotes closure),
\item $Cl(\bigcup U_i)=S^n$.
\end{enumerate}
Then we call the complement $S^n-\bigcup U_i$ an $(n-1)$-dimensional
{\it Sierpinski curve} (abbreviated to $\mathcal S$-curve).
\end{Def}

\begin{Thm}[Cannon, 1973]
Let $X,Y$ be an arbitrary pair of $(n-1)$-dimensional $\mathcal
S$-curves ($n\neq 4$). Then we have:
\begin{itemize}
\item $X$ is homeomorphic to $Y$,
\item if $i:X\rightarrow S^n$ is an arbitrary embedding, then $i(X)\subset S^n$ is an
$(n-1)$-dimensional $\mathcal S$-curve.
\item if $h:X\rightarrow Y$ is an arbitrary homeomorphism,
then $h$ extends to a homeomorphism of the ambient $n$-dimensional
spheres.
\end{itemize}
\end{Thm}

\vskip 5pt

The scheme of the proof is now clear: considering the double $DM_i$
of the manifold $M_i$ across its boundary, we can view $\tilde M_i$
as a totally geodesic subset of $\widetilde{DM}_i$, and hence $\pinf
\tilde M_i$ as an embedded subset of $\pinf \widetilde{DM}_i \cong
S^{n-1}$.  If we can establish that $\pinf \tilde M_i$ is an
$(n-2)$-dimensional $\mathcal S$-curve, Cannon's theorem will
immediately imply that $\pinf \tilde M_1$ is homeomorphic to $\pinf
\tilde M_2$.  We now proceed to verify the four conditions of an
$(n-2)$-dimensional $\mathcal S$-curve for $\pinf \tilde M \subset
\pinf \widetilde{DM} \cong S^{n-1}$.

Let us first fix some notation: the collection $\{U_i\}$ will be the
connected components of $\pinf \widetilde{DM} - \pinf \tilde M$
inside $\pinf \widetilde{DM} \cong S^{n-1}$.  We will denote by
$\{Y_i\}$ the connected components of $\widetilde{DM} -\tilde M$.
Note that each $Cl(Y_i)$ intersects $\tilde M$ along a boundary
component, which is a totally geodesic codimension one submanifold
of $\widetilde{DM}$.  We will denote by $Z_i \subset \partial \tilde
M$ the boundary component corresponding to $Y_i\subset
\widetilde{DM}-\tilde M$. Finally, we observe that each $U_i$ can be
identified with a corresponding $\pinf Y_i -\pinf Z_i$, for some
suitable component $Y_i$.

\vskip 5pt

\noindent {\bf Claim 1:}  The collection $\{U_i\}$ forms a null sequence.

\begin{Prf}
At the cost of rescaling the metric on $DM$, we may assume that the
sectional curvature is bounded above by $-1$, and hence that
$\widetilde{DM}$ is a $CAT(-1)$ space.  In this situation, Bourdon
\cite{B} defined a metric on $\pinf \widetilde{DM}$ inducing the
standard topology on $\pinf \widetilde{DM} \cong S^{n-1}$.  The
metric is given by:
$$d_\infty (p,q)=e^{-d(*,\gamma_{pq})}$$
where $\gamma_{pq}$ is the unique geodesic joining the points $p,q
\in \pinf \widetilde{DM}$, $*\in DM$ a chosen basepoint (and $d$
denotes the distance inside $\widetilde{DM}$).  Note that different
choices of basepoints result in metrics which are Lipschitz
equivalent.  For convenience, we will pick the basepoint $*$ to lie
in the interior of the lift $\tilde M$.

Now consider one of the components $U_i$, and let us try to estimate
$diam(U_i)$.  Note that given any two points $p,q \in Cl(U_i)$, we
have that the geodesic $\gamma_{pq} \subset Cl(Y_i)$, where $Y_i$ is
the component corresponding to $U_i$.  In particular, we see that
$d(*, \gamma_{pq})\geq d(*, Z_i)$, and hence that for any $p,q\in
Cl(U_i)$ we have the upper bound:
$$d_\infty (p,q)=e^{-d(*,\gamma_{pq})}\leq e^{-d(*, Z_i)}$$
Since $diam(U_i)$ is the supremum of $d_\infty(p,q)$, where $p,q \in
Cl(U_i)$, the above bound yields $diam(U_i)\leq e^{-d(*, Z_i)}$.  On
the other hand, since $\tilde M$ is the universal cover of a {\it
compact} negatively curved manifold with non-empty boundary, we have
that $\lim\{d(*, Z_i)\}= \infty$, where $Z_i$ ranges over the
boundary components of $\tilde M$.  This implies that the collection
$\{U_i\}$ forms a null sequence in $\pinf \widetilde{DM} \cong
S^{n-1}$, as desired.
\end{Prf}

\vskip 5pt

\noindent {\bf Claim 2:}  $S^{n-1}-U_i$ is an $(n-1)$-cell for each $i$.

\begin{Prf}
Recall that there exists a homeomorphism $\pi_x : S^{n-1}\cong \pinf
\widetilde{DM} \rightarrow T^1_x DM \cong S^{n-1}$, obtained by
mapping a point $p\in \pinf \widetilde{DM}$ to the unit vector
$\dot{\gamma}_{xp}(0)$, where $\gamma_{xp}$ is the unit speed
geodesic ray originating from $x$, in the direction $p\in \pinf
\widetilde{DM}$. Now let $U_i$ be given, and pick $x$ to lie on the
corresponding $Z_i$.  Note that under the homeomorphism $\pi_x$, we
have that $\pinf Z_i$ maps homeomorphically to a totally geodesic
$S^{n-2}\subset S^{n-1} \cong T^1_x \widetilde{DM}$, while the
subset $U_i$ maps homeomorphically to one of the open hemispheres
determined by $\pi_x(\pinf Z_i)$.  In particular, we see that $\pinf
\widetilde{DM} - U_i$ maps homeomorphically to one of the {\it
closed} hemispheres determined by $\pi_x(\pinf Z_i)$, and hence must
be an $(n-1)$-cell, as desired.
\end{Prf}

\vskip 5pt

\noindent {\bf Claim 3:}  $Cl(U_i)\cap Cl(U_j) =\emptyset$ for all $i\neq j$.

\begin{Prf}
Note that by definition we have that $U_i\cap U_j = \emptyset$, and
that $Cl(U_i)=U_i\cup \pinf Z_i$, $Cl(U_j)=U_j\cup \pinf Z_j$.
Hence it is sufficient to show that $\pinf Z_i\cap \pinf Z_j
=\emptyset$ for $i\neq j$ (since these are codimension one spheres
in $S^{n-1}\cong \pinf \widetilde {DM}$, with the $U_i$, $U_j$
connected components of the respective complements).  But a pair of
distinct boundary components of $\tilde M$, the universal cover of a
{\it compact} negatively curved manifold with non-empty totally
geodesic boundary, must diverge exponentially (with growth rate
bounded below in terms of the upper bound on sectional curvature).
In particular, no geodesic ray in $Z_i$ is within bounded Hausdorff
distance of a geodesic ray in $Z_j$, and hence the boundaries at
infinity are pairwise disjoint, as desired.
\end{Prf}

\vskip 5pt

\noindent {\bf Claim 4:}  $Cl(\bigcup U_i)=S^{n-1}$.

\begin{Prf}
Fix a point $x\in \tilde M$, and consider the homeomorphism $\pi_x:
S^{n-1}\cong \pinf \widetilde {DM} \rightarrow T^1_x \widetilde
{DM}\cong S^{n-1}$.  We will show that every point in $T^1_x
\widetilde {DM}\cong S^{n-1}$ can be approximated by a sequence of
points in $\pi_x(U_i)$.  This will imply that $T^1_x \widetilde {DM}
= Cl(\bigcup \pi_x(U_i))$, and since $\pi_x$ is a homeomorphism,
this will immediately imply Claim 4.

Now if $p\in T^1_x \widetilde {DM}$ lies in one of the $\pi_x(U_i)$,
we are done, so let us assume that $p\in T^1_x \widetilde {DM} -
\bigcup \pi_x(U_i)$. Let $\gamma$ be a unit speed geodesic ray
originating from $x$ with tangent vector $p$ at the point $x$.  Note
that we have that $\gamma \subset \tilde M \subset \widetilde {DM}$,
since we are assuming $p\in T^1_x \widetilde {DM} - \bigcup
\pi_x(U_i)$.  Now observe that $\tilde M$ is the universal cover of
a {\it compact} negatively curved manifold with non-empty totally
geodesic boundary, and hence there exists a constant $K$ with the
property that every point in $\tilde M$ is within distance $K$ of
$\partial \tilde M = \bigcup Z_i$ (for instance take $K=diam (M)$).

So for each integer $k\in \mathbb N$, we can find a point $y_k\in
\partial \tilde M$ satisfying $d(\gamma(k), y_k)\leq k$. Now observe
that if $\eta_k$ is the geodesic ray originating from $x$ and
passing through $y_k$, we have that $\eta_k(\infty) \in U_{i_k}$,
where $Z_{i_k}$ is the component of $\partial \tilde M$ containing
the point $y_k$.  This implies that $\dot {\eta} _k(0) \in
T^1_x\widetilde {DM}$ lies in the corresponding $\pi_x(U_{i_k}$,
i.e. that the sequence of vectors $\{\dot {\eta} _k(0)\} \subset
T^1_x\widetilde {DM}$ lies in the set $\bigcup \pi_x(U_i)$.  We are
left with establishing that $\lim \{\dot {\eta} _k(0)\} =p$.  To see
this, we need to estimate the angle between the geodesics $\eta_k$
and the geodesic $\gamma$.  But this is easy to do: consider the
geodesic triangle with vertices $(x, \gamma(k), y_k)$, and note that
$d(x,\gamma (k)=k$, while $d(\gamma(k), y_k)\leq K$.  Applying the
Alexandrov-Toponogov triangle comparison theorem, we see that the
angle $\angle (\dot{\eta}_k(0), \dot{\gamma}(0))$ is bounded above
by the angle of a comparison triangle in $\mathbb H^2$ (recall that
we assumed the metrics have been scaled to have upper bound $-1$ on
the sectional curvature).  But an easy calculation in hyperbolic
geometry shows that if one has a sequence of triangles in $\mathbb
H^2$ of the form $(A_k,B_k,C_k)$ with the property that
$d(A_k,B_k)=k$ and $d(B_k,C_k)\leq K$, then the angle at the vertex
$A_k$ tends to zero as $k$ tends to infinity.  This implies that
$\lim \{\angle (\dot{\eta}_k(0), \dot{\gamma}(0))\} =0$, and hence
completes the proof of Claim 4.
\end{Prf}

Appealing to Cannon's theorem, the four Claims above immediately
yield property (1) from Theorem 1.2: if $M_1, M_2$ are a pair of
compact $n$-dimensional negatively curved manifolds with non-empty,
totally geodesic boundary, then $\pinf \tilde M_1$ is homeomorphic
to $\pinf \tilde M_2$.  We now proceed to establish property (2):
under the hypotheses above, $\tilde M_1$ is homeomorphic to $\tilde
M_2$.  In order to do this, we pick a pair of points $p_i$ in the
{\it interior} of the respective $\tilde M_i$, and define subspaces
$C_i\subset \tilde M_i$ to be the union of all geodesic rays,
emanating from the respective $p_i$ to points in the corresponding
$\pinf \tilde M_i$.  Note that each $C_i$ is homeomorphic to the
{\it open cone} over $\pinf \tilde M_i$, that is to say the space
$\pinf \tilde M_i \times [0, \infty) /\sim$, where the equivalence
relation $\sim$ collapses $\pinf \tilde M_i \times \{0\}$ to a
point.  From property (1), we conclude that $C_1$ is homeomorphic to
$C_2$.  We now proceed to extend the homeomorphism between the
subsets $C_i\subset \tilde M_i$ to a homeomorphism between the
respective $\tilde M_i$.  We will denote by $C(Y)$ the open cone
over any topological space $Y$.

In order to extend the homeomorphism, let us view $\tilde M$ as a
subset in $\widetilde {DM}$.  Since each $\pinf Z_i \subset \pinf
\widetilde {DM} \cong S^{n-1}$ separates, the subset $C(\pinf
Z_i)\subset X$ separates $\widetilde {DM}$.  Let us denote by $H_i$
the unique connected component of $\pinf \widetilde {DM} - (Z_i \cup
C(\pinf Z_i))$ that contains both $Z_i$ and $C(\pinf Z_i)$ in its
closure.  Observe that $H_i\subset \tilde M$, and that we have a
decomposition of $\tilde M$ as $X \cup_i H_i$, where each $H_i$
attaches to $X$ along the boundary component $C(\pinf Z_i)\cong
\mathbb R ^{n-1}$.  Property (2) will now follow from the following:

\begin{Lem}
Each $H_i$ is homeomorphic to $[0,1] \times \mathbb R ^{n-1}$.
\end{Lem}

\begin{Prf}
Consider the point $x\in \widetilde {DM}$ from which we cone to
obtain $X= C(\pinf \tilde M)$, and observe that in $T_x^1\widetilde
{DM}\cong S^{n-1}$, we have that the set of directions to points in
$\pinf Z_i$ form an embedded codimension one sphere $S^{n-2}$ inside
$T_x^1\widetilde {DM}$.  Denoting by $S\subset T_x^1\widetilde{DM}$
this embedded codimension one sphere, we further observe that the
geodesics joining $x$ to any point in $H_i-C(\pinf Z_i)$ have the
property that they all lie in a common component $D$ of
$T_x^1\widetilde {DM} - S$ (and $D$ is homeomorphic to an open
$(n-1)$-dimensional cell).

Next we note that given any direction $v\in D$, the unit speed
geodesic ray $\gamma_v(t)$ has the property that its distance from
the subset $C(\pinf Z_i)$ tends to infinity as $t\rightarrow
\infty$.  Since the subset $H_i$ lies within finite Hausdorff
distance of $C(\pinf Z_i)$, this implies that the subset $R_v:=\{t
\hskip 3pt | \hskip 3pt \gamma _v(t) \in H_i\}$ is bounded.
Continuity ensures that the subset $R_v$ is closed inside
$[0,\infty)$, and it is clear that it is open at any $t\in R_v$ with
the property that $\gamma_v(t)\in H_i - Z_i$.  We now claim that the
set $R_v \subset [0, \infty)$ has only one boundary point.  Indeed,
if it had two such points $t_1<t_2$, then from the comments above,
we must have that $\gamma_v(t_1), \gamma_v(t_2)\in Z_i$ which
implies, since $Z_i$ is totally geodesic, that $\gamma_v\subset
Z_i$.  But we know that $\gamma_v(0)=x\notin Z_i$, yielding a
contradiction.

So we see that for each $v\in D$, the subset $R_v\subset [0,\infty)$
is a compact subset containing precisely one boundary point.  This
implies that it is a subinterval of $[0,\infty)$ of the type $[0,
\phi(v)]$, where $\phi(v)$ is a real number depending on the chosen
direction $v\in D$.  Now note that the function $\phi: D\rightarrow
[0, \infty)$ is a continuous function, tending to infinity as we
approach $\partial D = S$.  Furthermore, for each point $y\in
H_i-C(\pinf Z_i)$, there is unique $v(y)\in D$ and a unique $t_y\in
[0, \phi(v)]$ with the property that $\gamma _{v(y)}(t_y)=y$.

Now fix a homeomorphism $\rho$ from $D$ to the upper hemisphere in
the standard $(n-1)$-dimensional sphere $S^{(n-1)}\subset \mathbb
R^n$.  Construct a map $\bar \rho: H_i-C(\pinf Z_i)\rightarrow
\mathbb R^n$ by setting $\bar \rho(y):= \rho(y)\cdot \Phi_{\phi(v)}(t_y)$,
where the functions $\Phi_s:[0, s]\rightarrow [0,1]$ are homeomorphisms
varying continuously with respect to $s$, and having the property that
the map $\Phi_\infty: [0,\infty)\rightarrow [0,1)$ defined by $\Phi_\infty (t):=
\lim _{s\rightarrow \infty} \Phi_s(t)$ is a homeomorphism.  Observe that
the map $\bar \rho$ is a homeomorphism from $H_i-C(\pinf Z_i)$ to the
subset
$$\{ \vec{v}=(v_1, \ldots v_n)\in \mR ^n \hskip 3pt | \hskip 3pt v_n> 0,
||\vec{v}|| \leq 1\}$$
The homeomorphism $\bar \rho$ aligns (using $\rho$)
the directions $D$ pointing from
$x$ into the subspace $H_i-C(\pinf Z_i)$ with directions at the origin
in $\mR ^n$ pointing into the upper hemisphere, and then scales
(using the functions $\Phi_s$) the intervals so that each geodesic
segment $\gamma_v\cap [H_i-C(\pinf Z_i)]$ maps to the unit length
radial geodesic segment in the direction $\rho (v)$.  Now observe that
the choice of the scaling functions $\Phi_s$ implies that this
homeomorphism $\bar \rho$ {\it extends} to a homeomorphism
from $H_i$ to the subset:
$$\{ \vec{v}\in \mR ^n \hskip 3pt | \hskip 3pt v_n> 0,
||\vec{v}|| \leq 1\} \cup \{\vec{v}\in \mR ^n \hskip 3pt | \hskip 3pt v_n= 0,
||\vec{v}|| < 1\} $$
But the subset of $\mR^n$ described above is clearly homeomorphic to
$[0,1]\times \mR^{n-1}$, concluding the proof of Lemma 2.1.
\end{Prf}

To conclude the proof of Property (2), we take the homeomorphism
from $X_1$ to $X_2$.  Note that each connected component of $\tilde
M_1 - X_1$ is given by some $H_i\cong [0,1]\times \mathbb R^{n-1}$,
attached to a corresponding $C(\pinf Z_i)\subset X_1$, and
furthermore, the homeomorphism $\pinf \tilde M_1\rightarrow \pinf
\tilde M_2$ takes each $\pinf Z_i\subset \pinf \tilde M_1$
homeomorphically to a corresponding $\pinf Z_i^\prime \subset \tilde
M_2$.  This yields homeomorphisms between each $C(\pinf Z_i)\subset
X_1$ and the corresponding $C(\pinf Z_i^\prime)\subset X_2$.  On the
level of the corresponding $H_i\subset \tilde M_1$ and $H_i^\prime
\subset \tilde M_2$, this gives a homeomorphism between the subsets
$\{1\} \times \mathbb R^{n-1}$ in the respective topological
splittings $H_i \cong [0,1]\times \mathbb R^{n-1} \cong H_i^\prime$.
Extending in the obvious manner, we obtain {\it compatible}
homeomorphisms between the various $H_i\subset \tilde M_1$ and the
(bijectively associated) $H_i^\prime \subset \tilde M_2$.
Compatibility ensures that when we glue the $H_i$ to $X_1$, we still
obtain a homeomorphism onto the space obtained by gluing the
$H_i^\prime$ to $X_2$.  But the resulting spaces are $\tilde M_1$
and $\tilde M_2$ respectively, completing the proof of Property (2),
and hence of Theorem 1.2.

\begin{Rmk}
(1) As mentioned in the introduction, there is no hope of the
present argument giving a smooth classification of the universal
cover, as the homeomorphism is ``built'' from a homeomorphism
between the (non-manifold) boundary at infinity of the universal
cover.

(2) The hypothesis of {\it strict} negative curvature, rather than
non-positive curvature, is used in two places.  First of all, in the
proof of Claim 3, to ensure that distinct connected lifts of
boundary components yield subsets of the boundary at infinity that
are {\it disjoint}. In the case of non-positive curvature, there is
the possibility of two such distinct lifts of boundary components
containing geodesic rays that are asymptotic.  This can only occur
if there exists a (semi-infinite) flat strip isometric to
$[0,\infty)\times [0,r]$ (for some positive real number $r$) with
$[0,\infty) \times \{0\}$ mapping to one boundary component, and
$[0,\infty) \times \{r\}$ mapping to the other boundary component.
Hence to obtain Claim 3, one can weaken the curvature hypothesis
somewhat, by allowing non-positive curvature, but requiring the fact
that there do not exist any such flat strips. More problematic is
the use of strict negative curvature in the proof of Claim 1: in the
presence of zero curvature, we cannot use the Bourdon metric on the
boundary at infinity.  We can alternatively use the homeomorphic
projection to the unit tangent space at a point $x$ in the interior
of the manifold, but we still have a problem: in the zero curvature
setting, we can have a sequence of boundary components $Z_i$ with
$d(x,Z_i)\rightarrow \infty$, but with the $Z_i$ projecting to
subsets with diameter uniformly bounded away from zero.  It is not
clear what hypothesis would be needed to avoid this difficulty.

(3)  Note that compactness of the manifold $M^n$ is used
superficially in the proofs of Claims 1, 3 and 4.  Indeed the same
argument classifies topologically any simply connected Riemannian
$n$-dimensional ($n\neq 5$) manifold $X$ having the following
properties:
\begin{itemize}
\item $X$ has each boundary component totally geodesic and complete,
used to show Claim (2),
\item $X$ is semi-geodesically complete, i.e. every geodesic segment
with endpoints not lying on distinct boundary components is
extendible,
\item for any two distinct boundary components $Z_i,Z_j$, there exists
a constant $\epsilon_{ij}>0$ such that for all $p_i\in Z_i$, $p_j\in
Z_j$, we have the lower bound $d(p_i,p_j)\geq \epsilon_{ij}$, used in
Claim (3),
\item sectional curvature bounded above by $-a^2$ and bounded below
by $-b^2$, used in Claims (1) and (3),
\item the topological dimension of the boundary at infinity satisfies
$dim(\pinf X)<n-1$, used in Claim (4).
\end{itemize}
The first two bullets
are used to construct the complete simply connected Riemannian manifold
$\bar X \supset X$ by repeated reflections in the boundary components
(which serves as a substitute for $\widetilde {DM}$), and hence an inclusion
$\pinf X \subset \pinf \bar X \cong S^{n-1}$.
We leave it to the interested reader to verify that, with the conditions
mentioned above, the proofs of Claims (1)-(4) go through with
minimal changes (the author includes in the list above the number of the
Claims whose proofs require the specified bullet).  

\end{Rmk}

\section{Topological rigidity and applications.}

In this section, we explain the proofs of Theorem 1.1 and 1.2 from
the introduction. We first start with a definition:

\begin{Def}
Define the 1-{\it tripod} $T$ to be the topological space obtained
by taking the join of a one point set with a three point set. Denote
by $*$ the point in $T$ corresponding to the one point set. We
define the $n$-{\it tripod} ($n\geq 2$) to be the space $T \times
\mD^{n-1}$, and call the subset $*\times \mathbb D^{n-1}$ the {\it
spine} of the tripod $T\times \mathbb D^{n-1}$.  The subset $*\times
\mathbb D^{n-1}$ separates $T\times \mathbb D^{n-1}$ into three open
sets, which we call the {\it open leaves} of the tripod.  The union
of an open leaf with the spine will be called a {\it closed leaf} of
the tripod. We say that a point $p$ in a topological space $X$ is
{\it $n$-branching} provided there is a topological embedding
$f:T\times \mathbb D^{n-1} \longrightarrow X$ such that $p\in
f(*\times \mathbb D^{n-1}_\circ)$.
\end{Def}

It is clear that the property of being $n$-branching is invariant
under homeomorphisms.  Note that, in a simple, thick P-manifold of
dimension $n$, points in the codimension one strata are
automatically $n$-branching.  One can ask whether this property can
be detected at the level of the boundary at infinity.  This is the
content of the following:

\begin{Prop}[Characterization of branching points]
Let $X$ be an $n$-dimensional, simple, thick, negatively curved
P-manifold, and $p\in \pinf \tX$.  Then $p$ is $(n-1)$-branching if
and only if there exists a geodesic ray $\gamma$, entirely contained
in the lift of the branching locus, and satisfying $\gamma
(\infty)=p$.
\end{Prop}

\begin{Prf}
First observe that if $p\in \pinf \tX$ coincides with $\gamma
(\infty)$, for some $\gamma$ entirely contained in a connected
component $\mathcal B$ of the lift of the branching locus, then from
the thickness hypothesis, there exist $\geq 3$ lifts of chambers
that contain $\gamma$ in their common intersection $\mathcal B$.
Focusing on three such lifts of chambers, call them $Y_1,Y^\p_1,
Y^{\p\p}_1$, we can successively extend each of these in the
following manner: form subspaces $Y_{i+1}, Y^\p_{i+1},
Y^{\p\p}_{i+1}$ from the subspaces $Y_i,Y^\p_i, Y^{\p\p}_i$ by
choosing, for each boundary component of $Y_i,Y^\p_i, Y^{\p\p}_i$
{\it distinct} from $\mathcal B$, an incident lift of a chamber
(note that each boundary component is a connected component of the
lift of the branching locus).  Finally, set $Y_\infty := \cup_i
Y_i$, and similarly for $Y^\p_\infty, Y^{\p\p}_\infty$.  Now observe
that, by construction, the three subsets $Y_{\infty}, Y^\p_{\infty},
Y^{\p\p}_{\infty}$ have the following properties:
\begin{itemize}
\item they are totally geodesic subsets of $\tX$,
\item their pairwise intersection is precisely $\mathcal B$, their
(common, totally geodesic) boundary component,
\item doubling them across their boundary $\mathcal B$ results in a
simply connected, negatively curved,
complete Riemmanian manifold.
\end{itemize}
The first property ensures that the boundary at infinity of the
space $Y_{\infty}\cup Y^\p_{\infty}\cup Y^{\p\p}_{\infty}$ embeds in
$\pinf \tX$.  The third property ensures that $\pinf Y_{\infty}
\cong \pinf Y^\p_{\infty} \cong \pinf Y^{\p\p}_{\infty}\cong \mathbb
D ^{n-1}$.  The second property ensures that $S^{n-2}\cong \pinf
\mathcal B\subset \pinf \tX$ coincides with the boundary of the
three embedded $\mathbb D^{n-1}$.  Since $p\in \pinf \mathcal B$,
this immediately implies that $p$ is $(n-1)$-branching, yielding one
of the two desired implications.

Conversely, assume that $p \in \tX$ is {\it not} of the form $\gamma
(\infty)$, where $\gamma$ is contained entirely in a connected
component of the lift of the branching locus.  Consider a geodesic
ray $\gamma$ satisfying $\gamma (\infty)=p$, and note that there are
two possibilities:
\begin{itemize}
\item there exists a connected lift $W$ of a chamber with the property
that $\gamma$ eventually lies in the {\it interior} of
$W$, and is {\it not} asymptotic to any boundary component of $W$, or
\item $\gamma$ intersects infinitely many connected lifts of chambers.
\end{itemize}
In both these cases, we would like to argue that $p$ {\it cannot} be
$(n-1)$-branching.

Let us consider the first of these two cases, and assume that there
exists an embedding $f: T\times \mD ^{n-2}\rightarrow \pinf \tX$
satisfying $p\in f(\{*\}\times \mD_\circ ^{n-2})$.  Picking a point
$x$ in the interior of $W$, one can consider the composition
$\pi_x\circ f: T\times \mD^{n-2}\rightarrow lk_x \cong S^{n-1}$,
where $lk_x$ denotes a small enough $\epsilon$-ball centered at the
point $x$, and the map $\pi_x$ is induced by geodesic retraction.
Note that the map $\pi_x$ is {\it not} injective: the points in
$lk_x$ where $\pi_x$ is injective coincides with $\pi_x(\pinf W)$
(i.e. for every $q\in \pinf W$, we have $\pi_x^{-1}(\pi_xq) =
\{q\}$, and the latter are the only points in $\pinf \tX$ with this
property).  Note that, from Theorem 1.2, along with part (2) of
Cannon's theorem, this subset of injective points $I\subset lk_x$ is
an $(n-2)$-dimensional Sierpinski curve.  Furthermore, the
hypothesis on the point $p$ ensures that $\pi_xp$ does {\it not} lie
on one of the boundary spheres of the $(n-2)$-dimensional Sierpinski
curve $I$.  But now in \cite{L2} the following result was
established:

\vskip 5pt

\noindent {\bf Theorem:}  Let $F: T\times \mathbb D^{n-2}
\rightarrow S^{n-1}$ be a continuous map, and assume that the sphere
$S^{n-1}$ contains an $(n-2)$-dimensional Sierpinski curve $I$.  Let
$\{U_i\}$ be the collection of embedded open $(n-1)$-cells whose
complement yield $I$, and let $Inj(F)\subset S^{n-1}$ denote the
subset of points in the target where the map $F$ is injective. Then
$F(\{*\}\times \mD_\circ ^{n-2}) \cap [I-\cup_i (\partial U_i)] \neq
\emptyset$, implies that $[\cup_i(\partial U_i)] -Inj(F) \neq
\emptyset$.  In other words, this forces the existence of a point in
some $\partial U_i$ which has {\it at least two} pre-images under
$F$.

\vskip 5pt

Actually, in \cite{L2} this Theorem was proved using {\it purely
topological arguments} under some further hypotheses on the open
cells $U_i$.  But parts (1) and (3) of Cannon's Theorem allows the exact same
proof to apply in the more general setting, just by composing with a
homeomorphism taking the arbitrary Sierpinski curve to the one used
in the proof in \cite{L2}.

To conclude, we apply the Theorem above to the composite map
$F:=\pi_x\circ f: T\times \mathbb D^{n-2}\rightarrow lk_x$. The
point $f^{-1}(p)\in \{*\}\times \mD_\circ^{n-2}$ has image lying in
$I-\cup _i(\partial U_i)$, which tells us that $F(\{*\}\times
\mD_\circ ^{n-2}) \cap [I-\cup_i (\partial U_i)] \neq \emptyset$.
The Theorem implies that there exists a point $q$ in some $\partial
U_i\subset I$ which has {\it at least two} pre-images under the
composite map $F=\pi_x\circ f$.  Since the map $\pi_x$ is actually
{\it injective} on the set $I$, this implies that the map $f$ had to
have two pre-images at the point $\pi_x^{-1}(q)\in \pinf \tX$,
contradicting the fact that $f$ was an embedding.  This resolves the
first of the two possible cases.

For the second of the two cases (where the geodesic ray $\gamma$
passes through infinitely many lifts of chambers), a simple
separation argument (see Sections 3.2, 3.3 in \cite{L2}) shows that
if there exists a branching point of the second type, there must
also exist a branching point of the first type.  But we saw above
that there cannot exist any branching points of the first type. This
concludes the proof of Proposition 3.1.
\end{Prf}

Now given the characterization of branching points, let us see how
to show Theorem 1.1.  So assume that we are given a pair $X_1,X_2$
of simple, thick, negatively curved P-manifolds of dimension $n\geq
6$, and that we are told that $\pi_1(X_1)\cong \pi_1(X_2)$.  This
immediately implies that $\tX_1$ is quasi-isometric to $\tX_2$, and
hence that $\pinf \tX_1$ is {\it homeomorphic} to $\pinf \tX_2$.
Let $\mathcal B_i$ denote the union, in each respective $\pinf
\tX_i$, of the boundaries at infinity of the individual connected
components of the lift of the branching locus.  Note that each
$\mathcal B_i$ is a union of countably many, pairwise disjoint,
embedded $S^{n-2}$ inside $\pinf \tX_i$ (each $S^{n-2}$ arising as
the boundary at infinity of a single connected component of the lift
of the branching locus) .  Now the characterization of branching
points in Proposition 3.1 implies that, under the homeomorphism
between $\pinf \tX_1$ and $\pinf \tX_2$, we have that $\mathcal B_1$
must map homeomorphically to $\mathcal B_2$.

In particular, connected components of $\mathcal B_1$ must map
homeomorphically to connected components of $\mathcal B_2$.  A
result of Sierpinski \cite{Si} implies that the connected components
in each case are precisely the individual $S^{n-2}$ in the countable
union.  This yields a bijection between connected components of the
lift of the branching locus in the respective $\tX_i$.  Furthermore,
the homeomorphism must restrict to a homeomorphism between the {\it
complements} of the $\mathcal B_i$ in the respective $\pinf \tX_i$.
The connected components of this complement are either:
\begin{itemize}
\item isolated points, corresponding to $\gamma(\infty)$, where $\gamma$
is a geodesic ray passing through infinitely many connected lifts of chambers, and
\item components with $\geq 2$ points, which are in bijective correspondance
with connected lifts of chambers in the respective $\tX_i$ (see \cite[Section 3.2]{L2}).
\end{itemize}
This yields a bijective correspondence between lifts of chambers in $\tX_1$ and
lifts of chambers in $\tX_2$.  Furthermore, the {\it closure} of the components containing
$\geq 2$ correspond canonically with $\pinf W_i$, where $W_i$ is the bijectively
associated connected lift of a chamber.

Now recall that the homeomorphisms
between $\pinf \tX_1$ and $\pinf \tX_2$ has the additional property that
it is {\it equivariant} with respect to the respective $\pi_1(X_i)$
actions on the $\pinf \tX_i$.  We also have the following Lemma relating
the dynamics on $\pinf \tX$ with the action on $\tX$ (the argument is identical
to that given in \cite[pg. 212]{L1}) :

\begin{Lem}
Let $B_i$ be a connected component of the lift of the branching locus in
$\tX$, and let $W_i$ be a connected lift of a chamber in $\tX$.  Then we have:
\begin{itemize}
\item $Stab_{\pi_1(X)}(B_i) = Stab_{\pi_1(X)}(\pinf B_i)$, and
\item $Stab_{\pi_1(X)}(W_i) = Stab_{\pi_1(X)}(\pinf W_i)$.
\end{itemize}
where the action on the left hand side is the obvious action of $\pi_1(X)$ on $\tX$
by deck transformations, and the action on the right hand side is the induced action
of $\pi_1(X)$ on $\pinf \tX$.
\end{Lem}

Observe that equivariance of the homeomorphism implies that the bijective
correspondence between connected lifts of chambers descends to a bijective
correspondence between the chambers in $X_1$ and the chambers in $X_2$
(since two connected lifts of chambers cover the same chamber in $X_i$
if and only if the two lifts have stabilizers which are conjugate in $\pi_1(X_i)$).
Similarly, the bijective correspondance between connected components of
the lifts of the branching loci descends to a bijective correspondence between
the connected components of the branching loci in $X_1$ with those in $X_2$.
Furthermore, by equivariance of the homeomorphism, we have that chambers
(or connected components of the branching loci) that are bijectively identified
have isomorphic fundamental groups.  Separation arguments identical to the
ones in \cite[Lemmas 2.1-2.4]{L1} ensures that the bijective correspondence
also preserves the incidence relation between chambers and components of
the codimension one strata (and that the isomorphisms between the
various fundamental groups respect the incidence structure).

To conclude, we apply the celebrated Farrell-Jones topological
rigidity theorem for non-positively curved manifolds \cite{FJ}. This
implies that, corresponding to the bijections between chambers (and
components of the branching loci), one has {\it homeomorphisms}
between the corresponding chambers that induce the isomorphisms on
the level of the fundamental groups.  Note that, a priori, the
various homeomorphisms between chambers might not be compatible with
the gluing maps.  But by construction, the attaching maps all induce
the same maps on the fundamental group $\pi_1(B_i)$ of each
individual component $B_i$ of the branching locus.  By
Farrell-Jones, this implies that the restriction to $B_i$ of the
maps induced by the various homeomorphisms of incident chambers are
{\it all pairwise isotopic}.  Hence at the cost of deforming the
homeomorphism in a collared neighborhood of the boundary of each
chamber, we may assume that the homeomorphisms respect the gluing
maps.  But attaching together these individual homeomorphisms on
chambers now induces a {\it globally defined} homeomorphism from
$X_1$ to $X_2$.  This concludes the sketch of Theorem 1.1.


To obtain Corollary 1.1, we merely note that the generalized
Seifert-Van Kampen theorem implies that both $\pi_1(X_i)$ can be
expressed as the direct limit of a diagram of groups, with vertex
groups given by the fundamental groups of the chambers (and of the
components of the branching locus), and edge morphisms induced by
the inclusion of the components of the branching locus into the
incident chambers. Now an abstract isomorphism between the direct
limits corresponds to an isomorphism from $\pi_1(X_1)$ to
$\pi_1(X_2)$.  From Theorem 1.1, this isomorphism is induced by a
homeomorphism from $X_1$ to $X_2$, and hence must take chambers to
chambers and components of the branching locus to components of the
branching locus. This implies the existence of isomorphism between
the groups attached to the vertices in the diagram for $\pi_1(X_1)$
to the groups attached to the corresponding vertices in the diagram
for $\pi_1(X_2)$.  Furthermore, these isomorphisms commute (up to
inner automorphisms, due to choice of base points) with the
corresponding edge morphisms.  But this is precisely the definition
of diagram rigidity.  This concludes the sketch of Corollary 1.1.

Next let us explain the argument for Corollary 1.2.  Since the
space $X$ is a $K(\Gamma,1)$, any injection $i:\Gamma
\hookrightarrow \Gamma$ with image of finite index yields a
finte cover $\hat i: \bar X\rightarrow
X$ with $\pi_1(\bar X)\cong \Gamma$, and $\hat i (\pi_1(\bar
X))=i(\Gamma)$.  Now Theorem 1.2 implies that $\bar X$ is
homeomorphic to $X$, so this yields a covering map $\hat i:
X\rightarrow X$, whose degree coincides with the index of the group
$i(\Gamma)$ in $\Gamma$.  Hence it is sufficient to show that this
covering has degree one.  But we know that $X$ contains a chamber
with a non-zero characteristic number.  Since there are finitely
many chambers, consider a chamber $W$ for which this characteristic
number has the largest possible magnitude $|r|\neq 0$.  Then we know
that under a covering of degree $d$, characteristic numbers scale by
the degree, so we conclude that the pre-image chamber $\hat
i^{-1}(W)$ has characteristic number of magnitude $d\cdot |r|$.  By
maximality of the characteristic number of $W$, we conclude that
$d=1$, as desired.  Note that in this argument, it is crucial that the
image $i(\Gamma)$ has finite index in $\Gamma$.  If this is not the
case, then the covering space $\bar X$ is non-compact.  Since
compactness was an essential ingredient in the proof of Theorem 1.1,
one cannot conclude in this situation that $\bar X$ is homeomorphic to $X$.

Finally, for Corollary 1.3, take any element $\alpha\in Out(\Gamma)$.  Then
there exists an element $\bar \alpha \in Aut(\Gamma)$ which projects to $\alpha$
under the canonical map $Aut(\Gamma)\twoheadrightarrow Out(\Gamma)$.
From Theorem 1.1, we have a self-homeomorphism $\phi\in Homeo(X)$
with the property that $\phi_*=\alpha$, concluding the proof of Corollary 1.3.


\begin{thebibliography}{000bitem}




\bibitem[B]{B} M. Bourdon,
`Structure conforme au bord et flot geodesiques d'un
$CAT(-1)$-espace.' {\it Enseign. Math. (2)} 41 (1995), pp. 63--102.


\bibitem[BH]{BH} M.R. Bridson and A. Haefliger,
{\it Metric spaces of non-positive curvature} (Springer-Verlag,
Berlin, 1999).

\bibitem[C]{C} J.W. Cannon, `A positional characterization of the $(n-1)$-dimensional Sierpi\'nski curve in $S\sp{n}(n\not=4)$',
{\it Fund. Math.} 79 (1973), pp. 107--112.

\bibitem[FJ]{FJ} F.T. Farrell and L.E. Jones, `A topological analogue of Mostow's rigidity
theorem', {\it J. Amer. Math. Soc.} 2 (1989), pp. 257--370.


\bibitem[L1]{L1} J.-F. Lafont, `Rigidity result for certain 3-dimensional singular spaces
and their fundamental groups', {\it Geom. Dedicata} 109 (2004), pp.
197--219.

\bibitem[L2]{L2} J.-F. Lafont, `Strong Jordan separation and applications to rigidity',
to appear in {\it J. London Math. Soc.}

\bibitem[L3]{L3} J.-F. Lafont, `Rigidity of geometric amalgamations of free
groups', to appear in {\it J. Pure App. Alg.}

\bibitem[M]{M} G.D. Mostow, {\it Strong rigidity of locally symmetric spaces} (Princeton University Press,
Princeton, N.J., 1973).

\bibitem[Se]{Se} Z. Sela, `Structure and rigidity in (Gromov) hyperbolic groups and
discrete groups in rank $1$ Lie groups. II.', {\it Geom. Funct. Anal.} 7 (1997), pp. 561-593.

\bibitem[Si]{Si} W. Sierpinski, `Un th\'eor\`eme sur les continus', {\it Tohoku Math. Journ.} 13 (1918),
pp. 300-303.

\bibitem[W]{W} G.T. Whyburn, `Topological characterization of the Sierpi\'nski curve', {\it
Fund. Math.} 45 (1958), pp. 320-324.


\end{thebibliography}
\end{document}